\documentclass{amsart}

\usepackage{amscd}
\usepackage{amssymb}

\newcommand{\nr}{\mathrm{nr}}

\newcommand{\Zm}[1]{\Z/{#1}}

\newcommand{\injects}{\hookrightarrow}
\newcommand{\ra}{\rightarrow}

\newcommand{\norm}{\mathrm{norm}}

\DeclareMathOperator{\tr}{Tr} 
\DeclareMathOperator{\Id}{Id} 
\DeclareMathOperator{\Gal}{Gal}

\DeclareMathOperator{\aut}{Aut}

\newcommand{\QZ}{(\Q/\Z)'}

\newcommand{\Htnr}{H^3_{\mathrm{nr}}(BG)_{\mathrm{norm}}}
\newcommand{\Hdnr}{H^d_\nr}

\newcommand{\QZt}{\QZ(2)}
\newcommand{\basemu}{\boldsymbol{\mu}}

\newcommand{\mmu}[1]{\basemu_{#1}}
\newcommand{\mmud}[1]{\basemu_{#1}^{\otimes 2}}

\newcommand{\Fx}{F^\ast}
\newcommand{\Lx}{L^\ast}

\newcommand{\qform}[1]{{\langle{#1}\rangle}}                   

\DeclareMathOperator{\chr}{char} 
\DeclareMathOperator{\res}{res} \DeclareMathOperator{\cores}{cor}
\DeclareMathOperator{\Spin}{Spin} \DeclareMathOperator{\End}{End}

\newcommand{\s}{\sigma}

\renewcommand{\mmu}[1]{\boldsymbol{\mu}_{#1}}

\newcommand{\ut}{{\underline{t}}}

\newcommand{\GLt}{\widetilde{GL}}
\newcommand{\ot}{\otimes}

\newcommand{\oD}{{^1\hskip-.2em D_4}}

\newcommand{\iiiD}{{^3\hskip-.2em D_4}}
\newcommand{\viD}{{^6\hskip-.2em D_4}}

\newcommand{\oE}{{^1\hskip-.2em E_6}}
\newcommand{\dE}{{^2\hskip-.2em E_6}}

\newcommand{\e}{\varepsilon}

\newcommand{\jordmat}[6]{\left( \begin{array}{ccc} #1&#6&\cdot \\ \cdot&#2&#4 \\ #5&\cdot&#3 \end{array} \right) }
\newcommand{\sjordmat}[6]{\left( \begin{smallmatrix} #1&#6&\cdot \\ \cdot&#2&#4 \\ #5&\cdot&#3 \end{smallmatrix} \right) }

\newcommand{\sbasjord}{\sjordmat{\e_0}{\e_1}{\e_2}{c_0}{c_1}{c_2}}

\newcommand{\trio}[1]{#1^{\times\negthinspace 3}}
\newcommand{\trion}[1]{(#1_0, #1_1, #1_2)}

\newcommand{\diag}{\mathrm{diag}\,}

\DeclareMathOperator{\Inv}{Inv}
\newcommand{\Invt}{\Inv^3}
\newcommand{\Invtnr}{\Inv^3_{\nr}}

\newcommand{\C}{\mathfrak{C}}
\newcommand{\Z}{\mathbb{Z}}
\newcommand{\Q}{\mathbb{Q}}

\newcommand{\n}{\mathfrak{n}}

\newcommand{\la}{\lambda}

\newtheorem{thm}[equation]{Theorem}
\newtheorem{lem}[equation]{Lemma}

\newtheorem{prop}[equation]{Proposition}

\newtheorem{mainthm}[equation]{Main Theorem}

\newtheorem{app}[equation]{Application}

\newtheorem{silem}[equation]{Strongly Inner Lemma}

\theoremstyle{definition}

\newtheorem{eg}[equation]{Example}

\theoremstyle{remark}

\newtheorem{rmk}[equation]{Remark}

\numberwithin{equation}{section}

\makeatletter
\newenvironment{borel}[1]%
{\smallskip \noindent\refstepcounter{equation}{\bf \theequation.}{\bf{ #1.}}}%
{\smallskip \global\@ignoretrue} \makeatother

\makeatletter
\newenvironment{borel*}%
{\smallskip \noindent\refstepcounter{equation}{\bf \theequation.}}%
{\smallskip \global\@ignoretrue} \makeatother

\newcommand{\borelless}{%
\noindent\refstepcounter{equation}\textbf{\theequation.}~}

\begin{document}

\title[Unramified cohomology]%
{Unramified cohomology of classifying varieties for exceptional
simply connected groups}
\author{R.~Skip Garibaldi}
\date{18 April 2003}
\address{Department of Mathematics \& Computer Science, Emory University, Atlanta, GA 30322}
\email{skip@member.ams.org}
\urladdr{http://www.mathcs.emory.edu/{\textasciitilde}skip/}

\begin{abstract}
Let $BG$ be a classifying variety for an exceptional simple simply connected
algebraic group $G$.  We compute the degree 3 unramified Galois
cohomology of $BG$ with values in $\QZt$ over an arbitrary
field $F$. Combined with a paper by Merkurjev, this completes the
computation of these cohomology groups for $G$ semisimple simply
connected over all fields.

These computations provide another example of a simple simply connected group $G$ such that $BG$ is not stably rational.
\end{abstract}

\maketitle

\bigskip

Let $G$ be an algebraic group over a field $F$ with an embedding
$\rho \!: G \injects SL_n$ over $F$.  We call $BG$ a {\em
classifying space} of $G$.  We will compute the {\em unramified cohomology} of $BG$, defined
as follows.  We write $H^d(F)$ for the Galois cohomology group $H^d(\Gal(F), \QZ(d-1))$,
where $\QZ(d-1) = \varinjlim_n
\mmu{n}^{\otimes (d-1)}$ for $n$ not divisible by the characteristic
of $F$. For each $d
\ge 2$, define $H^d_\nr(BG/F)$ (or simply $\Hdnr(BG)$) to be the intersection of the kernels
of the residue homomorphisms
\begin{equation} \label{res.def}
\partial_v \!: H^d(F(BG)) \ra H^{d-1}(F(v))
\end{equation}
as $v$ ranges over the discrete valuations of $F(BG)$ over $F$.
The natural homomorphism $H^d(F) \ra \Hdnr(BG/F)$ is
split by evaluation at the distinguished point of $BG$; this gives a direct sum decomposition of $H^d_\nr(BG)$, and we denote the complement of $H^d(F)$ by $\Hdnr(BG)_\norm$.  This group depends only on $G$ and $F$ \cite[2.3]{M:ur}.

The purpose of this paper is to complete the computation of
$H^3_\nr(BG/F)_\norm$ for $G$ semisimple simply connected and $F$
 arbitrary.  The computation is quickly reduced to the case where $G$ is
simple simply connected \cite[\S 4]{M:ur}.  In \cite{M:ur},
$H^3_\nr(BG)_\norm$ was computed for $G$ simple and classical.  We
compute it for the remaining cases, where $G$ is {\em
exceptional}, that is, where $G$ is of type $G_2$, $\iiiD$,
$\viD$, $F_4$, $E_6$, $E_7$, or $E_8$.

\begin{mainthm} \label{MT}
Suppose that $G$ is a simple simply connected {\em exceptional}
algebraic group defined over a field $F$. Then
\[
H^3_\nr(BG)_\norm = \begin{cases} \Z/2 & \begin{cases}\parbox{2.5in}{if $\chr F \ne 2$, $G$ is of type
$\iiiD$, and $G$ has a nontrivial
Tits algebra} \end{cases} \\
0 & \text{otherwise.}
\end{cases}
\]
\end{mainthm}

(See \ref{char} for an explanation of the characteristic $\ne 2$ hypothesis.)


The general motivation for studying $\Hdnr(X)$ is that it can sometimes detect if $X$ is not stably rational, see \cite[pp.~35--39]{CT:birat}.  It was an open question whether $BG$ is stably rational for $G$ semisimple simply connected.
The first counterexamples were provided in \cite{M:ur}, where Merkurjev exhibited groups $G$ of type $^2\!A_n$, $^2\!D_3$, and $^1\!D_4$ with $\Htnr \ne 0$.
The results here give another such $G$, this time of type $\iiiD$.

\medskip

Our basic tool is that one can compute $\Htnr$ by inspecting the ramification of the Rost invariant of $G$, see \cite{M:ur} or \ref{ram.connect}.  Many questions are settled by hopping along the chain of inclusions
   \[
   G_2 \subset D_4 \subset F_4 \subset E_6 \subset E_7 \subset E_8
   \]
of split groups, see \S\ref{trivcase} and \ref{first.borel}. 

The most interesting part of the proof is where we show that the mod 4 portion of the Rost invariant is ramified for groups of type $\dE$.  We prove  (in \ref{D4E6}) that every isotropic trialitarian group embeds in a group of type $\dE$ with trivial Tits algebras.  This settles the question, since the mod 4 portion of the Rost invariant for groups of type $\viD$ is easily shown to be ramified \eqref{6D4}.  The proof of \ref{D4E6} uses Galois descent and  interpretations of exceptional groups as acting on nonassociative algebras.  

\begin{rmk}
Computations of $H^d_\nr(X/F)$ in the literature for $X$ a smooth
variety (e.g., a classifying space) typically assume that $F$ is
algebraically closed. The examples of nonrational classifying
varieties $BG$ provided here and in \cite{M:ur} require that $F$
is {\em not} algebraically closed.
\end{rmk}

\section{Vocabulary} \label{invsec} \label{vocab}

We say that an (affine) algebraic group is {\em simple} if it is $\ne 1$,
is connected, and has no nontrivial connected normal subgroups
over an algebraic closure. (These groups are often called
``absolutely almost simple''.)  Simple groups are classified in,
e.g., \cite[Ch.~VI]{KMRT}. We say that a group is {\em of type
$T_n$} if it is simple with root system of type $T_n$ and {\em of
type $^tT_n$} if additionally the absolute Galois
group of $F$ acts as a group of automorphisms of order $t$ on the
Dynkin diagram.

Let $V$ be a finite-dimensional irreducible representation of an algebraic group $G$ over $F$.
The $F$-algebra $\End_G(V)$ is  a skew field by Schur's
Lemma, and it is finite-dimensional over $F$; it is called a {\em
Tits algebra} for $G$.  If it is a (commutative) field, we say
that it is {\em trivial}.

The {\em Dynkin index} $n_G$ of a simple simply connected algebraic group $G$ is a natural number which depends only
upon the type of $G$ and the (Schur) indices of its Tits algebras.
We will repeatedly make us of the fact that the value of $n_G$ is
known and can be found by looking in \cite[App.~B]{MG} or \cite[pp.~437--442]{KMRT}.
The value of $n_G$ for $G$ exceptional is:
  \[
  \begin{array}{||c||c|c|c|c|c|c|c|c||} \hline
&    &^{3,6}\negthinspace D_4\text{, all Tits}          &^{3,6}\negthinspace D_4\text{, some Tits}& &&&&\\
\text{type of $G$}&G_2 &\text{alg's trivial}&\text{alg's nontrivial}&F_4&\oE&\dE&E_7&E_8\\
\hline n_G&2 & 6 & 12 & 6 & 6 & 12 & 12 & 60 \\ \hline
\end{array}
   \]
   
We have functors $H^1(*, G)$ and $H^3(*)$ which take a field extension of $F$ and give a pointed set and abelian group respectively.  A {\em degree 3 invariant of $G$ with values in $\QZ(2)$} is a morphism of functors
   \[   
  H^1(*, G) \longrightarrow H^3(*)
  \]
which takes base points to base points.  Such invariants are often called ``normalized''.  We write $\Invt(G)$ for the abelian group of such invariants.
(Clearly, this definition makes sense for every algebraic group $G$ over $F$.)
   
Let $p = \chr F$ if the characteristic is prime and $p = 1$ otherwise.  Write $n_G = p^k n'_G$, where $n'_G$ is a natural number prime to $p$.  
The group $\Invt(G)$ is cyclic of order $n'_G$.  It has a canonical generator $r_G$ which we call the {\em Rost invariant} of $G$.  (This is the prime-to-$p$ part of what is called the Rost invariant in \cite{MG}.  I do not know how to define a residue map for the $p$-primary part.)

For $\alpha \!: H \ra G$ a map between simple simply connected
algebraic groups over $F$, there is a positive integer $n_\alpha$
called the {\em Rost multiplier} or ``Dynkin index'' of $\alpha$,  see \cite[\S7]{MG}.  It has the properties: $n_H$ divides $n_\alpha n_G$ and 
for every extension $E$ of $F$, the composition
\[
\begin{CD}
H^1(E, H) @>{\alpha}>> H^1(E, G) @>{r_{G}}>> H^3(E)
\end{CD}
\]
is $n_\alpha r_{H}$.  

\section{$A_2 \subset D_4$} \label{A2D4.sec}

\borelless \label{A2D4.hyp}
In this section, we assume that $F$ contains a primitive cube root of unity and hence has characteristic $\ne 3$.  Let $L$ be a cubic Galois extension of $F$; by Kummer theory it is obtained by adjoining a cube root of some element $\la \in \Fx$.  We write $(\la)$ for the corresponding class in $\Fx / F^{*3} = H^1(F, \mmu3)$, where $\mmu3$ is the algebraic group of cube roots of unity.

The short exact sequence $1 \ra \mmu3 \ra SL_3 \ra PGL_3 \ra 1$ induces a connecting homomorphism
$\delta \!: H^1(F, PGL_3) \ra H^2(F, \mmu3)$.
   
\begin{lem} \label{A2D4} Continue the hypotheses of \ref{A2D4.hyp}.   Let
$G$ be the quasi-split simply connected group of type $\iiiD$ associated
with the extension $L/F$.  Then $G$ contains a subgroup isomorphic to $PGL_3$ such that 
for every extension $E$ of $F$ the diagram
\[
\begin{CD}
H^1(E, PGL_3) @>{x \mapsto \delta(x) \cup (\la)}>> H^3(E, \mmud{3}) \\
@VVV @VVV \\
H^1(E, G) @>{r_{G}}>> H^3(E)
\end{CD}
\]
commutes up to sign, where the arrow on the right comes from the natural map $\mmud{3} \ra \QZt$.
\end{lem}

\begin{proof}
We have maps
  \begin{equation} \label{A2D4.map}
   \begin{CD}
  PGL_3 \times \mmu3 @>>> \Spin_8 \rtimes \Zm3 @>>> F_4
   \end{CD}
  \end{equation}
where $F_4$ denotes the split algebraic group of that type and we identify $\mmu3$ with $\Zm3$ using the primitive cube root of unity in $F$.  The first map comes from the fact that $PGL_3 = \aut(M_3(F))$ preserves the subspace of trace zero elements in $M_3(F)$, see \cite[pp.~504, 505]{KMRT}.  The second map comes from the Springer decomposition of Albert algebras, see \cite[38.7]{KMRT}.  The group $\Zm3$ acts on $\Spin_8$ in a manner which cyclically permutes the vector and half-spin representations and fixes $PGL_3$ elementwise.  The map $\Spin_8 \ra F_4$ has Rost multiplier 1.

When we twist the groups in \eqref{A2D4.map} by $(\la)$ and restrict to connected components, we obtain the sequence
  \begin{equation} \label{A2D4.twist}
  \begin{CD}
  PGL_3 @>{\iota}>> G @>{\sigma}>> F_4
  \end{CD}
  \end{equation}
The map $H^1(\s \iota)$ sends a class $[A] \in H^1(E, PGL_3) $ of a central simple $E$-algebra of degree 3 to the class of the first Tits construction $[J(A, \la)] \in H^1(F, F_4)$, see \cite[39.9]{KMRT}.
The composition $r_{F_4} \circ H^1(\s \iota)$ is, up to sign, the composition of $\delta$ with the cup product $\cdot \cup (\la)$ by \cite[p.~537]{KMRT}.  Since $\s$ has Rost multiplier 1, $r_{F_4} \circ H^1(\s) = r_{G}$, and the lemma is proved.
\end{proof}

\section{$\iiiD, \viD \subset \dE$} \label{E6sec}

A simple algebraic group is said to be {\em trialitarian} if it is of type
$\iiiD$ or $\viD$.
In this section, we assume that $F$ has characteristic $\ne 2$;
our goal is to prove the following:

\begin{thm} \label{D4E6}  {\rm ($\chr F \ne 2$)}
Let $T$ be a
trialitarian simply connected group over $F$ which is
$F$-isotropic but not $F$-quasi-split.  Let $K$ be a quadratic extension of $F$ such that $T$ is $K$-quasi-split.  Then there exists a simply
connected group $G$ of type $\dE$ over $F$ such that
   \begin{enumerate}
   \item all of $G$'s Tits algebras are trivial;
   \item $G$ is of type $\oE$ over $K$; and
   \item $T$ is a subgroup of $G$ with Rost multiplier $1$.
   \end{enumerate}
\end{thm}

Given a $T$ as in the first sentence of \ref{D4E6}, such a $K$ always exists by \cite[0.1]{G:iso}. We postpone the proof of the theorem until the end of this section.

\begin{rmk}
If $T$ is simply connected trialitarian and $F$-quasi-split, one can do better.  In that case, $T$ is isomorphic to a subgroup of the split $F_4$ (see e.g.\ \cite[4.4]{G:iso}), hence $T$ is isomorphic to a subgroup of every simply connected and $F$-quasi-split group $G$ of type $E_6$.  Moreover, the inclusion $T \ra G$ has Rost multiplier 1 since over an algebraic closure of $F$ the inclusion arises from the natural inclusion of root systems $D_4 \subset E_6$.
\end{rmk}

\borelless \label{Cayley.basis} 
Let $\C$ be the split Cayley algebra over $F$ with canonical involution (a.k.a. ``conjugation'') $\pi_\C$.  Fix the basis $u_1$, $\ldots$, $u_8$ of $\C$ as in \cite{G:iso} and \cite{G:rinv} so that 
the bilinear norm form is $\n(x, y) = x \pi_\C(y) + y \pi_\C(x)$ is given by
$\n(u_i, u_j) = \delta_{(i+j),9}$ (Kronecker delta).
Write $\s$ for the involution on $GL(\C)$ which is adjoint for $\n$.
  
Let $R$ denote the subgroup of $GL(\C)^{\times 3}$ consisting of so-called {\em related triples} of proper similitudes of $\n$, see \cite[\S1]{G:iso} or \cite[\S35]{KMRT} for a definition.
This group is reductive with center of rank 2; its derived subgroup consists of triples $\ut = \trion{t}$ with $t_i \in SO(\n)$ for all $i$ and is isomorphic to $\Spin_8$ \cite[35.7]{KMRT}.

The group $S_3 = \qform{ r, \pi \mid r^3 = \pi^2 = 1, \pi r = r^2 \pi}$ acts on $R$ via
   \[
   {^r \ut} = (t_1, t_2, t_0) \quad \text{and} \quad 
    {^\pi (t_0, t_1, t_2)} = (\pi_\C t_0 \pi_\C, \pi_\C t_2 \pi_\C, \pi_\C t_1 \pi_\C).
    \]
Define $R \rtimes S_3$ to be the Cartesian product $R \times S_3$ with multiplication
   \[
   (\ut, \alpha) \cdot (\ut', \beta) = (\ut \cdot {^\alpha \ut'}, \alpha \beta).
   \]

The split Albert algebra $J$
has  underlying vector space the matrices in $M_3(\C)$ fixed
by the conjugate transpose.  With that in mind, we may write a
general element of $J$ as in \eqref{alb.gen} below where $\e_i \in F$, $c_i \in \C$, and the entries given as $\cdot$
are forced by symmetry.  The algebra $J$ has a canonically determined norm
form; write $\Inv(J)$ for the group of norm isometries.

There is an injection $g \!: R \rtimes S_3 \injects \Inv(J)$ defined by  
\begin{gather*}
g_\ut \sbasjord := \sjordmat{\mu(t_0)^{-1} \e_0}{\mu(t_1)^{-1}
\e_1}{\mu(t_2)^{-1} \e_2}{t_0(c_0)}{t_1(c_1)}{t_2(c_2)}, \\
g_r \sbasjord = 
\jordmat{\e_1}{\e_2}{\e_0}{c_1}{c_2}{c_0}, \quad
\text{and} \quad g_\pi \sbasjord =
\jordmat{\e_0}{\e_2}{\e_1}{\overline{c_0}}{\overline{c_2}}{\overline{c_1}}.
\end{gather*}

\borel{Construction of a quasi-split $\dE$}
The algebra $J$ is also endowed with a nondegenerate symmetric bilinear form $s$ defined by 
   \[
   s(x,y) = \tr_J(xy) = \sum_{i=0}^2 \left[ \e_i \nu_i + \n(c_i, d_i) \right]
   \]
for 
   \begin{equation} \label{alb.gen}
   x = \sbasjord, \qquad y = \sjordmat{\nu_0}{\nu_1}{\nu_2}{d_0}{d_1}{d_2}.
   \end{equation}
For each $f \in GL(J)$, there is a
unique $f^\dag \in GL(J)$ such that $s(f(x), f^\dag(y)) = s(x,y)$
for all $x, y \in J$. 

The map $f \mapsto f^\dag$ restricts to automorphisms of $\Inv(J)$ and $R \rtimes S_3$ defined over $F$.  We have
   \[
   r^\dag = r, \quad \pi^\dag = \pi, \quad \text{and} \quad \ut^\dag = (\s(t_0)^{-1}, \s(t_1)^{-1}, \s(t_2)^{-1}).
    \]

Fix a quadratic extension $K$ of $F$ with $\Gal(K/F) = \qform{\iota}$. 
We define the group
$E^K_6$ and $H$ to be the groups $\Inv(J)$ and $R$ with twisted $\iota$-actions: For $f$ a $K$-point, we set $^\iota\negthinspace f = \iota f^\dag \iota$
where the action on the left is the new action and
juxtaposition denotes the usual action.  The group $E^K_6$
is quasi-split of type $\dE$.

\medskip
\borelless {\sl Proof of Theorem \ref{D4E6}.}
Since $T$ is isotropic and not quasi-split, it  has a Tits algebra which is a nonsplit quaternion algebra $Q$ over a cubic extension $L$ of $F$ \cite[0.1]{G:iso}.  By \cite[43.9]{KMRT}, $Q$ is isomorphic to the
quaternion $L$-algebra $(a, b)_L$ for some $b \in \Fx$ and $a \in
\Lx$ such that $N_{L/F}(a) = 1$.  Since $Q$ is nonsplit, $b$ is not a square in $F$; we set $K = F(\sqrt{b})$.
Since $Q$ is not split over $L$, it is not split over $L^c$ by
\cite[3.2]{G:iso}.  In particular, $L^c$ does not contain a square
root of $b$, so $P = K \otimes_F L^c$ is a quadratic field extension of $L^c$.

To simplify our argument, we assume that $L$ is not Galois over
$F$, so $\Gal(L^c/F)$ is isomorphic to $S_3$.  (This is the case that will be used in the rest of the paper.  The other case --- where $L$
is Galois over $F$ --- is only easier.) Then, the group $\Gal(P/F)$ is
isomorphic to $S_3 \times \mmu2 \cong \Z/6 \rtimes \mmu2$, where
the factor of $\Z/6$ corresponds to the subgroup $\Gal(P/\Delta)$
for $\Delta$ the unique quadratic extension of $F$ in $L^c$. We
fix generators $\zeta := (r^{-1}, -1) \in S_3 \times \mmu2$ (which generates the
copy of $\Z/6$) and $\tau = (\pi, 1)$ (which generates the copy
of $\mmu2$ in $\Z/6 \rtimes \mmu2$ corresponding to
$\Gal(P/L(\sqrt{b}))$).

We construct the group $G$ by descent as follows.  The group $E^K_6$ --- and hence $H \rtimes S_3$ --- is a closed subgroup of $GL(V)$ for some $F$-vector space $V$.  (Specifically, $E^K_6$ is the group of algebra automorphisms of a Brown algebra with underlying vector space $V$, cf.\ \cite[2.9(2)]{G:struct}.)  We call an additive homomorphism $f \!: V \ot P \ra V \ot P$ {\em $\varpi$-semilinear} if there is some $\varpi \in \Gal(P/F)$ such that
   \[
   f(p v) = \varpi(p) f(v) \quad \text{for all $p \in P$ and $v \in V \ot P$.}
   \]
Let $\GLt(V)$ denote the (abstract) group of such maps $f$.  We define a group homomorphism $\phi \!: \Gal(P/F) \ra \GLt(V)$ such that $\phi(\varpi)$ is $\varpi$-semilinear for all $\varpi$ and $T(F)$ and $G(F)$ are the subgroups of $\Spin_8(P)$ and $E^K_6(P)$ commuting with $\phi(\varpi)$ for all $\varpi$.

Define  $\ut =
\trion{t} \in \trio{GL(\C)}$ by setting $t_i = m_i P$, for
\[
m_i = \diag(1, \rho^i(a), -\rho^i(a), \rho^{i+2}(a)^{-1},
\rho^{i+1}(a)^{-1}, -1, 1, \rho^i(a))
\]
with $\rho := \zeta^2$, and $P$ the matrix permuting the basis
vectors as $(1\, 2)(3\,6)(4\,5)(7\,8)$, for the basis of $\C$ fixed in \ref{Cayley.basis} above. 
Since $N_{L/F}(a) = 1$, $\ut$ is a related triple by \cite[1.5(3), 1.6, 1.8]{G:iso}.
Set
   \begin{equation} \label{basic}
   \phi(\zeta) = \ut r \zeta \quad \text{and} \quad \phi(\tau) = \pi \tau.
   \end{equation}
We have
  \[
  \zeta \ut  \zeta^{-1} = (\s(\zeta(t_0)), \s(\zeta(t_1)), \s(\zeta(t_2)))^{-1} = (\s(t_2), \s(t_0), \s(t_1))^{-1} = r^{-1} \ut^{-1} r,
 \]
since $\s(t_i) = t_i$ for all $i$.
Hence 
     \begin{gather}
  \phi(\zeta)^2 = r^2 \zeta^2, \quad 
  \phi(\zeta)^3 =\ut \zeta^3, \quad \text{and} \notag \\
  \phi(\zeta)^6 = \Id_{V \ot P}.\label{basic.zeta}
  \end{gather}
Since $\pi$ and $\tau$ commute, we have 
  \begin{equation} \label{basic.tau}
  \phi(\tau)^2 = \Id_{V \ot P}.
  \end{equation}
Since 
   \[
   \phi(\zeta)^5 = \phi(\zeta)^2 \cdot \phi(\zeta)^3 =  \ut r^2 \zeta^5,
    \]
and $\pi \tau$ and $\ut$ commute,
it is easy to verify that 
    \begin{equation} \label{basic.comm}
   \phi(\tau) \phi(\zeta) = \phi(\zeta)^5 \phi(\tau).
   \end{equation}

Equations \eqref{basic.zeta}, \eqref{basic.tau}, and \eqref{basic.comm} give that \eqref{basic} defines a homorphism $\phi \!: \Gal(P/F) \ra \GLt(V)$.  Then the set map $z \!: \Gal(P/F) \ra H \rtimes S_3$ defined by
    $z_\varpi := \phi(\varpi) \varpi^{-1}$
 is in fact a 1-cocycle.
(This correspondence between groups of semilinear transformations and 1-cocycles is well-explained in \cite[\S3]{Jac:tri}.)
Set $G$ to be the
twisted group $(E^K_6)_z$; it automatically satisfies (2).  Since $z$ takes values in the simply connected group $E^K_6$, (1) holds.

Since the values of $z$ normalize the subgroup $\Spin_8$ of $E^K_6$, the
twisted group $(\Spin_8)_z$ is a subgroup of $G$.  The inclusion has Rost multiplier 1 since the inclusion $\Spin_8 \subset E_6^K$ over an algebraic closure arises from the natural inclusion of root systems $D_4 \subset E_6$.
The restriction of $z$ to $\Spin_8$ is the descent given in
\cite[4.7]{G:iso} to construct $T$, i.e., $(\Spin_8)_z$ is
isomorphic to $T$, hence (3). $\hfill\qed$

\begin{rmk}
The isotropic group $G$ occurring in Theorem \ref{D4E6}
is typically not quasi-split, even over $L$.  This can be seen by
examining the mod 2 portion of the Rost invariant for $(z) \in H^1(P/L, E^K_6)$, which
is typically nontrivial by \cite[6.7]{G:rinv}.
\end{rmk}

\section{A construction} \label{constsec}\label{wacky.sec}
\newcommand{\Lhat}{\widehat{L}}

The purpose of this section is to construct a suitable extension of $F$ over which we may apply \ref{D4E6}:

\begin{app} {\rm ($\chr F \ne 2$)} \label{wacky.app}
Let $G$ be a quasi-split simply connected group of type $\dE$ over a field $F$.  There exists
   \begin{itemize}
   \item an extension $F'$ of $F$ such that $G$ is of type $\dE$ over $F'$,
   \item an isotropic, non-quasi-split group $T'$ over $F'$ of type $\viD$ 
      \item a strongly inner form $G'$ of $G$ over $F'$
   \end {itemize}
such that $T'$ is a subgroup of $G'$ with Rost multiplier $1$.
\end{app}

A group is a {\em strongly inner form} of the simply connected group $G$ if it is obtained from $G$ by twisting by a 1-cocycle in $Z^1(F, G)$.

The following arguments are
somewhat simpler than previously, thanks to suggestions by Adrian
Wadsworth.

\begin{lem} {\rm ($\chr F \ne 2$)}\label{wacky.1}
For $p, q \in \Fx$, the ring $L = F(t)[x] / (x^3 + px + qt)$ is a
separable cubic field extension of $F(t)$ which is not Galois over $F(t)$.
There is a prolongation of the $t$-adic valuation on $F(t)$ to $L$
which is unramified with residue degree $1$ and with respect to
which $x$ has value $1$.
\end{lem}

\begin{proof}
If $L$ is not a field, then there is some $a \in F(t)$ such that
$a^3 + pa + qt = 0$.  Since $a$ is integral over the UFD $F[t]$,
it belongs to $F[t]$, so it makes sense to speak of the degree of
$a$.  In particular, at least two of the terms $a^3$, $pa$, and
$qt$ must have the same degree, which is also the maximum of the
degrees.  This implies that $a$ cannot have positive degree.  But
then $qt$, with degree 1, is the unique term of maximal degree,
which is a contradiction.

Since $p$, $q$ are in $\Fx$, the discriminant $-4p^3 - 27q^2 t^2$
of $L$ is not 0, hence $L$ is separable over $F$.  An argument similar to the one in the preceding paragraph shows that the discriminant is not a square in $F(t)$: Any square root $b \in F(t)$ of the discriminant would belong to $F[t]$ and have degree 1.  Then
the coefficient of $t$ in $b^2$ would be nonzero.  Thus $L$ is not
Galois over $F(t)$.

Hensel's Lemma gives that $x^3 + px + tq$ has a linear factor of
the form $x - \pi$ in $F((t))[x]$, where $\pi$ has $t$-adic value 1. The
map $x \mapsto \pi$ gives an isomorphism of $L$ with the subfield
$F(t)(\pi)$ of $F((t))$, and the $t$-adic valuation obviously
extends to $L$ so that $x$ has value 1.  Since $F((t))$ is the
completion of $F(t)$ with respect to the $t$-adic valuation and
hence is unramified with residue degree 1, the claims about
ramification and residue degree of our prolongation to $L$ follow.
\end{proof}

\begin{lem} {\rm ($\chr F \ne 2$)}\label{wacky.2}
Let $p, b \in \Fx$ be such that the quaternion algebra $(p, b)_F$
is nonsplit. Let $L$ be as in Lemma \ref{wacky.1}.  Then the
quaternion algebra $(x, b)_L$ is nonsplit and is not isomorphic to $(-qt, b)_L$.
\end{lem}

\begin{proof}
Since $N_{L/F(t)}(x) = -qt$, the corestriction of $(x,b)_L$ down
to $F(t)$ is Brauer-equivalent to $(-qt, b)_{F(t)}$.  This algebra
is split if and only if the quadratic form $\qform{1, -b, qt}$ is
isotropic over $F(t)$.  Over the completion $F((t))$, this form
has residue forms $\qform{1, -b}$ and $\qform{q}$.  Since the
algebra $(p, b)_F$ is nonsplit, the first form is anisotropic,
hence $\qform{1, -b, qt}$ is anisotropic over $F((t))$ by
Springer's Theorem.  Thus $(-qt, b)_{F(t)}$ is
nonsplit, and hence so is $(x, b)_L$.

For the sake of contradiction, suppose that $(x,b)_L$ is
 isomorphic
to $(-qt, b)_L$, i.e., the algebra $(-xqt, b)_L$ is
split.  Since
\[
-x(qt) = -x(-x^3 - px) = x^4 + px^2 \equiv x^2 + p \mod{L^{\ast
2}},
\]
the algebra $(x^2 + p, b)_L$ is split.

Let $\Lhat$ be a completion of $L$ with respect to the
prolongation of the $t$-adic valuation on $F(t)$ given by Lemma
\ref{wacky.1}. The norm of $(x^2 + p, b)_L$ is the form $\qform{1,
-(x^2 + p), -b, b(x^2 + p)}$ over $L$.  Since $x$ has value 1,
over $\Lhat$ this form has one residue form $\qform{1, -p, -b,
bp}$ over the residue field $F$.  This is the norm of the algebra
$(p, b)_F$, which is anisotropic because the algebra is nonsplit.
By Springer's Theorem, the norm of $(x^2 + p, b)_L$ is anisotropic
over $\Lhat$, hence the algebra is not $L$-split, which
contradicts our assumption that $(x, b)_L$ is isomorphic to $(-qt, b)_L$.
\end{proof}

\borelless {\sl Proof of \ref{wacky.app}.}
Let $F(\sqrt{b})$ be the quadratic extension which splits the given quasi-split group $G$.
Let $F_0 := F(p, t)$ for $p$, $t$ indeterminates.   Set $L_0 := F_0(t)[x] / (x^3 + px + t)$; it is a separable cubic non-Galois extension of $F_0$ by \ref{wacky.1}.  Set $F'$ to be the function field of the Severi-Brauer variety of the quaternion algebra $(-t, b)_{F_0}$.  Since $F_0$ is algebraically closed in $F'$, $L' := L_0 \otimes_{F_0} F'$ is a field which is cubic and not Galois over $F'$; it is the function field of the Severi-Brauer variety of $(-t, b)_{L_0}$.

Since $b$ is not a square in $F$, the quaternion algebras $(p, b)_{F_0}$ is not split.  
By \ref{wacky.2},  $(x, b)_{L_0}$ is nonsplit and is not isomorphic to $(-t, b)_{L_0}$, hence $(x, b)_{L'} \cong (x, b)_{L_0} \otimes_{L_0} L'$ is not split  by a well-known theorem of Amitsur.  The corestriction $\cores_{L'/F'}(x,b)_{L'}$ is $(-t, b)_{F'}$, which is split.  Thus there is a simply connected isotropic trialitarian group $T'$ over $F'$ with nontrivial Tits
algebra $(x, b)_{L'}$ \cite[4.7]{G:iso}.  It is of type $\viD$ since $L'$ is not
Galois over $F'$.  By \ref{D4E6}, there is a simple simply connected
group $G'$ over $L'$ of type $\dE$ which is a strongly inner form
of $G_{L'}$ such that $T'$ is a subgroup of $G'$ with Rost multiplier 1.$\hfill\qed$

\section{Ramification} \label{ramsec}

Let $G$ be an algebraic group over $F$.
We say that an invariant $a \in \Invt(G)$ is {\em unramified} if the composition
  \[
\begin{CD}
H^1(E((t)), G) @>{a}>> H^3(E((t))) @>\partial>> H^2(E)
\end{CD}
\]
is 0 for every field extension $E$ of $F$.  Otherwise we say that $a$ is {\em ramified}.  The following example is typical:

\begin{eg}\label{A2D4.eg} ($\chr F \ne 3$)
Let $F$ be a field with a primitive cube root of unity, and let $L = F(\la^{1/3})$ be a cubic Galois extension of $F$.  We claim that the invariant $a \in \Inv^3(PGL_3)$ from \S\ref{A2D4.sec} given by $a(x) = \delta(x) \cup (\la)$ is ramified.  

The set $H^1(F, PGL_3)$ classifies degree 3 cyclic central simple algebras $(C, a)$ for $C$ a cubic Galois extension of $F$ and $a \in F^*$.  Let $[C] \in H^1(F, \Zm3)$ denote the class corresponding to $C$.  Then
$\delta(C, a) \cup (\la)$ is  $\pm[C] \cup (a) \cup (\la)$ in $H^3(F, \mmud{3})$.
Taking $E = F(u)$ for $u$ an indeterminate and $C = E((t))(t^{1/3})$ a cubic Galois extension of $E((t))$, we have
   \[
   \partial \left[ \delta(C, u) \cup (\la) \right] = \pm (u) \cup (\la) \quad \in H^2(E).
   \]
This is nonzero in $H^2(E)$ since $u$ is not a norm from the extension $E(\la^{1/3})/E$.  Hence $a$ is ramified, as claimed.$\hfill\qed$
\end{eg}

For the rest of the section, we assume that $G$ is simple and simply connected.  We write $\Invtnr(G)$ for the subset of unramified invariants in $\Invt(G)$.  It is a subgroup since $\partial$ is a group homomorphism.

\begin{lem} \cite{M:ur} \label{ram.connect}
$H^3_\nr(BG)_\norm \cong \Invtnr(G)$.$\hfill\qed$
\end{lem}

\noindent In particular, $H^3_\nr(BG)_\norm$ is
necessarily finite, see \S\ref{vocab}.
\smallskip

\borelless \label{char} 
If $F$ has prime characteristic $p$, then multiplication by $p$ is an isomorphism of $\QZt$.  Hence $H^3(F)$, $\Inv^3(G)$, and --- by \ref{ram.connect} --- $H^3_\nr(BG/F)_\norm$ have no nontrivial $p$-torsion.

This explains the hypothesis ``$\chr F \ne 2$" in the Main Theorem: a simply connected group $G$ of type $\iiiD$ with nontrivial Tits algebras ``should'' have $\Htnr = \Z/2$, but this is impossible in characteristic 2.

\begin{silem} \label{silem}
Let $G$ be a simple simply connected group over $F$, and fix $z
\in Z^1(F, G)$.  The canonical identification $\Invt(G) = \Invt(G_z)$ defined by $r_G \leftrightarrow r_{G_z}$ restricts to an identification $\Invtnr(G) = \Invtnr(G_z)$.
\end{silem}

\begin{proof}
Let $E$ be an extension of $F$, and consider the diagram
\[
\begin{CD}
H^1(E((t)), G) @>{mr_{G}}>> H^3(E((t))) @>{\partial}>> H^2(E)\\
@V{\tau_z}V{\cong}V @V{\cdot - mr_G(z)}VV @|\\
H^1(E((t)), G_z)@>{mr_{G_z}}>> H^3(E((t)))
@>{\partial}>> H^2(E),
\end{CD}
\]
where $\tau_z$ is the twisting isomorphism. The left box commutes by
\cite[p.~76, Lem.~7]{Gille:inv} or \cite[1.7]{MPT}.  The right box commutes because
$\partial$ is a group homomorphism and $\partial(r_G(z)) =
0$.

Hence $mr_G$ is ramified if and only if $mr_{G_z}$ is.  Since
$G$ and $G_z$ are strongly inner forms of each other, they
have $n_G = n_{G_z}$, hence $n'_G = n'_{G_z}$.  The claim follows.
\end{proof}

\borel{Functoriality (homomorphisms)} \label{first.borel}
Let $\alpha \!: H \ra G$ be a morphism of algebraic groups.  Then $\alpha$ induces a natural map
   \[
   \alpha^* \!: \Invt(G) \ra \Invt(H)
   \]
which restricts to a homomorphism
   \[
  \alpha^*_\nr \!: \Invtnr(G) \ra \Invtnr(H).
  \]

Now suppose that $H$ and $G$ are simple simply connected.
If the Rost multiplier of $\alpha$ is 1, then $n_H$ divides $n_G$, hence $n'_G / n'_H$ divides $n_G/n_H$.  (See \S\ref{vocab} for definitions.)  Also, $\alpha^*$ is a surjection with kernel of order $n'_{G}/ n'_{H}$.  Then we have: {\em If $\Invtnr(H)$ is trivial, then $\Invtnr(G)$ is $(n_{G}/n_{H})$-torsion.}

\borel{Functoriality (scalar extension)} \label{scalar.ext}
Let $K$ be an extension field of $F$, and write $G_K$ for $G \times_F K$.  The restriction homomorphism
   \[
  \res_{K/F} \!: \Invt(G/F) \ra \Invt(G/K)
   \]
is the natural surjection $\Zm{n'_G} \ra \Zm{n'_{G_K}}$; its kernel
is the $(n'_G/n'_{G_K})$-torsion in $\Invt(G/F)$.
It restricts to a homomorphism
  \[
  (\res_{K/F})_\nr \!: \Invtnr(G/F) \ra \Invtnr(G/K).
  \]
The kernel of this map is killed by $n'_G / n'_{G_K}$, hence by $n_G / n_{G_K}$.
We have: {\em If $\Invtnr(G/K)$ is trivial, then $\Invtnr(G/F)$ is $(n_G/n_{G_K})$-torsion.}

\section{The case where $G$ has trivial Tits algebras}\label{trivcase}

In this section, we prove:

\begin{prop} \label{trivalg}
Suppose that $G$ is a simple simply connected exceptional algebraic group
defined over a field $F$.  If $G$ has only trivial Tits algebras, then
$\Invtnr(G) = 0$.
\end{prop}

That is, the Main Theorem holds for groups with only trivial Tits algebras by \ref{ram.connect}.
In proving the proposition, we may assume that $G$ is quasi-split by the Strongly Inner Lemma \ref{silem}.

The proposition still holds if the hypothesis ``exceptional'' is dropped; the classical groups are treated in \cite{M:ur}.

\borel{Type $\iiiD$} \label{trisplit}
Let $L$ be a separable cubic extension of $F$ over which
the quasi-split group $G$ is of type $\oD$.  Since $G$ is quasi-split, $n_{G_L} = 2$ and $\Invtnr(G/L) = 0$ by \cite[8.5]{M:ur}.  Hence $\Invtnr(G/F)$ is 3-torsion by \ref{scalar.ext}.  By \ref{char} we may assume that $\chr F \ne 3$.

Let $F'$ be the extension obtained from $F$ by adjoining (if not already in $F$) 
a primitive cube root of unity.  Then $G$ is still of type $\iiiD$ over $F'$ and the invariants $2r_G$ and $4r_G$ are ramified over $F'$ by \ref{A2D4} and \ref{A2D4.eg}.  Since $2r_G$ and $4r_G$ are the only nontrivial 3-torsion elements of $\Invt(G/F)$, we have shown that $\Invtnr(G)$ is 0.

\begin{lem}\label{6D4} Let $G$ be simple simply
connected of type $^6\!D_4$ over $F$.  Then $\Invtnr(G) =
0$.
\end{lem}

\begin{proof}
Let $\Delta$ be the unique quadratic extension of $F$ over which $G$ is of type $\iiiD$.  Let $K$ be a generic quasi-splitting field for $G$ over $\Delta$ as in
\cite{KR}.   Then $G$ is quasi-split of type $\iiiD$ over $K$, hence $n_{G_K} = 6$ and $\Invtnr(G/K) = 0$ by
the $\iiiD$ case \eqref{trisplit}.  We have that $\Invtnr(G/F)$ is 2-torsion by \ref{scalar.ext}.

Let $L$ be a cubic extension of $F$ over which $G$ is of type
$^2\!D_4$.  Then $n_{G_L} =  2$ or $4$ and $\Invtnr(G/L) = 0$ by \cite[8.5]{M:ur}.  Since $n_G = 6$ or 12 (as $n_{G_L} = 2$ or 4), 
$\Invtnr(G/F)$ is 3-torsion.

Combining the two previous paragraphs, we find $\Invtnr(G/F) = 0$.
\end{proof}

\borel{Type $G_2$} \label{G2} Here $n_G = 2$, so we may assume that $\chr F \ne 2$.  The Rost invariant for the split group of type $G_2$ is given explicitly in \cite[p.~441]{KMRT}.  It is the Elman-Lam invariant for 3-Pfister quadratic forms, which is clearly ramified.

\borel{Type $F_4$} \label{F4}
The split $G_2$ is contained in our split group $G$ of type $F_4$ with Rost multiplier 1.
Since $n_{F_4} = 6$ and $n_{G_2} = 2$, the group $\Invtnr(G)$ is $3$-torsion by \ref{first.borel} and the $G_2$ case \eqref{G2}.
In characteristic $\ne 3$, the mod 3 part of the Rost invariant is described in \cite[3.2]{PR:el}, and it is clearly ramified.  So $\Inv^3_\nr(G) = 0$.

\begin{borel}{Type $\oE$} \label{oE.case}
The split group $G$ of type $E_6$ contains a subgroup which is split of type $F_4$ \cite[14.20, 14.24]{Sp:jord} with Rost multiplier 1 \cite[2.4]{G:rinv}, and we have $n_{F_4} = n_G = 6$.
Hence $\Invtnr(G) = 0$ by \ref{first.borel} and the $F_4$ case (\ref{F4}).
\end{borel}

\begin{borel}{Type $\dE$} \label{dE.case}
The group $G$ is split by a quadratic extension and $n_G = 12$.  By \ref{scalar.ext} and the $\oE$ case \eqref{oE.case}, $\Invtnr(G)$ is 2-torsion, hence we may assume that $\chr F \ne 2$.

Let $G'$, $L'$, and $T'$ be as constructed in
\ref{wacky.app}.  By Lemma \ref{6D4}, $\Invtnr(T'/{L'}) =
0$.  Since $T'$ is isotropic and non-quasi-split, it has nontrivial Tits algebras \cite[5.6]{G:iso} and $n_{T'} = 12$.  Since $n_{G'} = 12$ and the inclusion $T' \subset G'$ has Rost multiplier 1,  we have $\Invtnr(G'/L') = 0$.  By the
Strongly Inner Lemma, $\Invtnr(G/{L'}) = 0$.  Since $G$ is
 of type $^2\!E_6$ over $F$ and $L'$, we have $n_G =
n_{G_{L'}}$, hence $\Invtnr(G/F) = 0$ by 
\ref{scalar.ext}.
\end{borel}

\begin{borel}{Type $E_7$} \label{E7.case}
The natural inclusion of root systems gives a split simply connected subgroup of type $E_6$ inside the split group $G$ of type $E_7$.  Since $n_G = 12$ and $n_{E_6} = 6$, $\Invtnr(G)$ is 2-torsion by \ref{first.borel} and the $\oE$ case \eqref{oE.case}.  Hence we may assume that $\chr F \ne 2$.

Set $F' = F(x)$ and $K = F'(\sqrt{x})$. There is a quasi-split
simply connected group $E^K_6$ over $F$ of type $\dE$ associated
with the extension $K/F'$; it injects into $G_{F'}$ with
Rost multiplier 1 \cite[\S 3]{G:rinv}.  Since $n_{E^K_6} = n_{G_{F'}}
= 12$, $\Invtnr(G/{F'}) = 0$ by \ref{first.borel} and
the $\dE$ case \eqref{dE.case}.  Since $n_G = 12$, we have $\Invtnr(G/F) = 0$.
\end{borel}

\borel{Type $E_8$}  As in the previous cases, we may assume that our group $G$
of type $E_8$ is actually split.  The natural
inclusion of root systems gives an embedding of a split simply connected group of
type $E_7$ in $G$ with Rost multiplier 1, so $\Invtnr(G)$ is 5-torsion by \ref{first.borel} and the $E_7$ case \eqref{E7.case}.  In particular, we may assume that $\chr F \ne 5$ by \ref{char}.

Let $F'$ be the extension obtained from $F$ by adjoining two indeterminates and (if necessary) a primitive 5th root of unity.  There is an $F'$-central division algebra $D$ of dimension $5^2$, namely the symbol algebra determined by the two indeterminates.  Arguing in a manner similar to \cite[\S1]{Gille:E8}, one finds a strongly inner form $G'$ of $G$ over $F'$ and an injection $SL_1(D) \hookrightarrow G'$ with Rost multiplier 1.
Now $n_{SL_1(D)} = 5$ \cite[11.5]{MG} and $\Inv^3_\nr(SL_1(D)/F') = 0$ (as can be seen from the explicit formula for the Rost invariant in \cite[p.~138]{MG}), hence $\Invtnr(G'/F')$ is 12-torsion by \ref{first.borel}.  By the Strongly Inner Lemma, $\Invtnr(G/F')$ is 12-torsion.  Since the Dynkin index of $G$ is 60 over $F$ and $F'$, $\Invtnr(G/F)$ is 12-torsion by \ref{scalar.ext}.

Combining the two preceding paragraphs gives that $\Invtnr(G)$ is 0.  This completes the proof of Prop.~\ref{trivalg}. $\hfill\qed$

\section{Proof of the main theorem} \label{3D4}

Let $G$ be as in the Main Theorem \ref{MT}.  If $G$ has only trivial
Tits algebras (e.g., $G$ is of type $G_2$, $F_4$, or $E_8$), then
the Main Theorem holds for $G$ by Prop.~\ref{trivalg} and Lemma \ref{ram.connect}.

If $G$ is of type $E_6$ or $E_7$, then we pick a generic
quasi-splitting field $K$ of $G$ over $F$.  We have $n_{G_K} =
n_G$ and $\Invtnr(G/K) = 0$ (by \ref{trivalg}), hence
$\Invtnr(G/F) = 0$ by \ref{scalar.ext}.

If $G$ is of type $^6\!D_4$, the Main Theorem holds by Lemma
\ref{6D4}.  The remaining case is where $G$ is of type $^3\!D_4$
with nontrivial Tits algebras.  We have $n_G = 12$ and as in the
proof of \ref{6D4},  $\Invtnr(G)$ is 2-torsion.  Hence we may assume that $\chr F \ne 2$. 

The only nontrivial 2-torsion element of $\Invt(G)$ is $6r_G$, so
we will complete the proof of the Main
Theorem if we show that $6r_G$ is unramified.  That is, if we show
that for every extension $E$ of $F$, the composition
\begin{equation} \label{tricomp}
\begin{CD}
H^1(E((t)), G) @>{6r_{G_{E((t))}}}>> H^3(E((t))) @>{\partial}>> H^2(E)
\end{CD}
\end{equation}
is trivial.

If $G$ is of type $\oD$ over $E$, then $n_{G_E} = 2$ or 4 \cite[15.4]{MG}.  Hence $6r_{G_{E((t))}} = 2{r_{G_{E((t))}}}$ and the composition \eqref{tricomp} is trivial by \cite[8.2]{M:ur}.

Otherwise, $G_E$ is of type $\iiiD$.  That is, if $L$ is a cubic Galois extension of $F$ over which $G$ is of type $\oD$, the tensor product $L_E := L \otimes_F E$ is a cubic field extension of $E$.
We have a diagram
\[
\begin{CD}
H^1(E((t)), G) @>{6r_{G_{E((t))}}}>> H^3(E((t)))
@>\partial>> H^2(E) \\
@V{\res}VV @V{\res}VV @V{\res_{L_E/E}}VV \\
H^1(L_E((t)), G) @>{6r_{G_{L_E((t))}}}>> H^3(L_E((t))) @>\partial>> H^2(L_E).
\end{CD}
\]
The left box commutes because the Rost invariant is compatible
with restriction, and the right box commutes because the extension
$E((t)) \subset L_E((t))$ is unramified,
hence the whole diagram
commutes.  Fix a class $\alpha$ in $H^1(E((t)), G)$. 
Since $G_{L_E((t))}$ is of type $\oD$, the composition
of the two bottom arrows is 0 by the preceding paragraph, and the
image of $\alpha$ in $H^2(L_E)$ is 0.
Let $\beta \in H^2(E)$ be the image of $\alpha$; it is 2-torsion because $6r_G$ is 2-torsion.  Hence
   \[
   \beta = 3 \beta = \cores_{L_E/E} \res_{L_E/E}(\beta),
   \]
which is 0 by the commutativity of the diagram.
This shows that the composition
\eqref{tricomp} --- which is the top row of the diagram --- is 0 in this case.

Thus
$\Invtnr(G) = \Z/2$ for $G$ of type $^3\!D_4$ with nontrivial
Tits algebras when $\chr F \ne 2$.
This completes the proof of the Main Theorem \ref{MT}.
$\hfill\qed$

\subsection*{Acknowledgements} I thank Alexander Merkurjev for suggesting the problem and for many enlightening conversations.

\providecommand{\bysame}{\leavevmode\hbox to3em{\hrulefill}\thinspace}
\providecommand{\MR}{\relax\ifhmode\unskip\space\fi MR }
\providecommand{\MRhref}[2]{%
  \href{http://www.ams.org/mathscinet-getitem?mr=#1}{#2}
}
\providecommand{\href}[2]{#2}

\end{document}